\documentclass{amsart}
\usepackage{amsmath,amssymb,amsfonts}
\usepackage{eucal}
\usepackage{verbatim}

\usepackage{enumerate}

\numberwithin{equation}{section}

\theoremstyle{plain}
\newtheorem{theorem}[equation]{Theorem}

\newtheorem{cor}[equation]{Corollary}
\newtheorem{prop}[equation]{Proposition}

\theoremstyle{definition}
\newtheorem{definition}[equation]{Definition}

\newtheorem*{note}{Note}

\newtheorem*{thank}{Acknowledgments}

\input xy
\xyoption{all}

\newcommand{\Deltaop}{{\bf \Delta}^{op}}

\newcommand{\SSets}{\mathcal{SS}ets}

\newcommand{\Sets}{\mathcal Sets}

\newcommand{\Secat}{\mathcal Se \mathcal Cat}

\newcommand{\sk}{\text{sk}}

\newcommand{\SSpof}{\mathcal{SS} p_{\mathcal O,f}}
\newcommand{\SSpoc}{\mathcal{SS} p_{\mathcal O,c}}
\newcommand{\LSSpof}{\mathcal{LSS}p_{\mathcal O,f}}
\newcommand{\LSSpoc}{\mathcal{LSS}p_{\mathcal O,c}}

\begin{document}

\title{A Characterization of Fibrant Segal Categories}

\author[J.E. Bergner]{Julia E. Bergner}

\address{Kansas State University \\
138 Cardwell Hall \\
Manhattan, KS 66506}

\email{bergnerj@member.ams.org}

\date{\today}


\begin{abstract}
In this note we prove that Reedy fibrant Segal categories are
fibrant objects in the model category structure $\Secat_c$.
Combining this result with a previous one, we thus have that the
fibrant objects are precisely the Reedy fibrant Segal categories.
We also show that the analogous result holds for Segal categories
which are fibrant in the projective model structure on simplicial
spaces, considered as objects in the model structure $\Secat_f$.
\end{abstract}

\maketitle

\section{Introduction}

Segal categories are simplicial spaces which are meant to look
like simplicial categories, but their morphisms are only
associative up to (higher) homotopy. They were first described by
Dwyer, Kan and Smith, who called them special $\Deltaop$-spaces
\cite{dks}. They have since been generalized to the notion of
Segal $n$-categories, variants of which have been studied as
models for weak $n$-categories by Hirschowitz and Simpson
\cite{hs}, \cite{simpson} and Tamsamani \cite{tam}. Like
simplicial categories, they can be considered as models for
homotopy theories \cite[8.6]{thesis}.

In any model category structure, it is useful to have a precise
description of the fibrant and cofibrant objects, since these are
used to define the homotopy category.  Furthermore, various
constructions, such as that of mapping spaces, are only homotopy
invariant when the objects involved are both fibrant and
cofibrant. In the first of the model structures which we will
consider, $\Secat_c$, all objects are cofibrant, but a description
of the fibrant objects has not been so clear. In this note, we
give a complete characterization of them.


We begin by giving some precise definitions. Let ${\bf \Delta}$
denote the cosimplicial category, or category whose objects are
finite ordered sets $[n] = (0 \rightarrow 1 \rightarrow \ldots
 \rightarrow n)$ for $n \geq 0$ and whose morphisms are order-preserving maps
between them. Then $\Deltaop$ is the opposite of this category and
is called the simplicial category. Recall that a \emph{simplicial
set} $X$ is a functor $\Deltaop \rightarrow \Sets$. We will denote
the category of simplicial sets by $\SSets$. (In the course of
this paper we will sometimes refer to simplicial sets as spaces,
due to their homotopy-theoretic similarity with topological spaces
\cite[3.6.7]{hovey}.) We denote by $\Delta [n]$ the $n$-simplex
for each $n \geq 0$, by $\dot \Delta [n]$ its boundary, and by
$V[n,k]$ the boundary with the $k$th face removed. More details
about simplicial sets can be found in \cite[I]{gj}. A simplicial
set $X$ is \emph{discrete} if all elements of $X_n$ are degenerate
for $n>0$.  We denote by $|X|$ the topological space given by
geometric realization of the simplicial set $X$ \cite[I.2]{gj}.

More generally, recall that a \emph{simplicial object} in a
category $\mathcal C$ is a functor $\Deltaop \rightarrow \mathcal
C$.  In particular, a functor $\Deltaop \rightarrow \SSets$ is a
\emph{simplicial space} or \emph{bisimplicial set} \cite[IV]{gj}.
Given a simplicial set $X$, we denote by $X^t$ the simplicial
space such that $(X^t)_0$ is the discrete simplicial set $X_0$.
Also, given a simplicial space $W$, we denote by $\sk_n W$ the
$n$-skeleton of $W$, or the simplicial space generated by the
simplices of $W$ of degree less than or equal to $n$ \cite[\S
1]{reedy}. In this paper, our primary concern is the case where
the simplicial set in degree zero is discrete (or constant).


\begin{definition}
A \emph{Segal precategory} is a simplicial space $X$ such that
$X_0$ is a discrete simplicial set.
\end{definition}

Now note that for any simplicial space $X$ there is a \emph{Segal
map}
\[ \varphi_k: X_k \rightarrow \underbrace{X_1 \times_{X_0} \ldots
\times_{X_0} X_1}_k \] for each $k \geq 2$, which we define as
follows.  Let $\alpha^i:[1] \rightarrow [n]$ be the map in ${\bf
\Delta}$ such that $\alpha^i(0)=i$ and $\alpha^i(1)=i+1$, defined
for each $0 \leq i \leq n-1$. We can then define the dual maps
$\alpha_i:[n]\rightarrow [1]$ in $\Deltaop$.  Given a simplicial
space $X$, for each $k \geq 2$ the Segal map is defined to be the
map
\[ \varphi_k: X_k \rightarrow \underbrace{X_1 \times_{X_0} \cdots \times_{X_0} X_1}_k \]
induced by the maps
\[ X(\alpha_i):X_k \rightarrow X_1. \]

\begin{definition} \cite[\S 2]{hs} \label{SeCat}
A \emph{Segal category} $X$ is a Segal precategory such that the
Segal map $\varphi_k$ is a weak equivalence of simplicial sets for
each $k \geq 2$.
\end{definition}

For a Segal (pre)category $X$, we will frequently refer to the
discrete simplicial set $X_0$ as the set of ``objects" of $X$.

In \cite[5.1, 7.1]{thesis}, we show that there exist two different
model category structures on the category of all Segal
precategories in which the fibrant objects are Segal categories.
(We will discuss model category structures in more detail in the
next section.) For the first of these structures, $\Secat_c$, we
show in \cite[5.13]{thesis} that the fibrant objects are Segal
categories which are fibrant in the Reedy model category structure
on the category of all simplicial spaces. Here, we complete the
result and show that the converse holds as well, namely, that all
Reedy fibrant Segal categories are fibrant in $\Secat_c$.
Similarly, in the second model structure, $\Secat_f$, we show that
the fibrant objects are precisely the Segal categories which are
fibrant in the projective model category structure on simplicial
spaces.

It should be noted that these model category structures on Segal
precategories fit into a chain of Quillen equivalences between
various model structures.  The two are Quillen equivalent to one
another, as well as to a model structure on the category of
simplicial categories and to Rezk's complete Segal space model
structure on simplicial spaces, which Joyal and Tierney prove to
be Quillen equivalent to Joyal's model structure on
quasi-categories \cite{joyal}, \cite{jt}.  In doing so, they
actually obtain an alternate proof of Theorem \ref{main}. The
author's interest in comparing these model structures arose from
finding models for the homotopy theory of homotopy theories, a
project begun by Rezk \cite{rezk}.

\begin{thank}
I would like to thank Bertrand To\"{e}n for pointing out an error
in a previous proof of the main result, and the referee for
suggestions which led to an improved proof of Proposition
\ref{oc}.
\end{thank}

\section{Review of model category structures}

In this section, we give a brief review of model category
structures.  In particular, we discuss the Reedy model category
structure on the category of all simplicial spaces and the model
category structure $\Secat_c$ on the category of Segal
precategories.

Recall that a model category structure on a category $\mathcal C$
is a choice of three distinguished classes of morphisms:
fibrations, cofibrations, and weak equivalences. A (co)fibration
which is also a weak equivalence is an \emph{acyclic
(co)fibration}. With this choice of three classes of morphisms,
$\mathcal C$ is required to satisfy five axioms MC1-MC5 which can
be found in \cite[3.3]{ds}.  An object $X$ in a model category is
\emph{fibrant} if the unique map $X \rightarrow \ast$ to the
terminal object is a fibration. Dually, $X$ is \emph{cofibrant} if
the unique map from the initial
object $\phi \rightarrow X$ is a cofibration.  


There is a model category structure on the category of simplicial
sets in which the weak equivalences are the maps $f:X \rightarrow
Y$ such that the induced map $|f|:|X| \rightarrow |Y|$ of
topological spaces is a weak homotopy equivalence
\cite[I.11.3]{gj}.  We can then use this model category structure
to define model structures on the category of simplicial spaces.

A natural choice for the weak equivalences in the category of all
simplicial spaces is the class of levelwise weak equivalences of
simplicial sets. Namely, given two simplicial spaces $X$ and $Y$,
we define a map $f:X \rightarrow Y$ to be a weak equivalence if
and only if for each $n \geq 0$, the map $f_n:X_n \rightarrow Y_n$
is a weak equivalence of simplicial sets.

In the Reedy model category structure on simplicial spaces
\cite{reedy}, the weak equivalences are the levelwise weak
equivalences of simplicial sets.  We will denote the Reedy model
structure by $\SSets^{\Deltaop}_c$. While it is defined somewhat
differently, the Reedy model category structure on simplicial
spaces is exactly the same as the injective model category
structure on this same category, in which the cofibrations are
defined to be levelwise cofibrations of simplicial sets
\cite[15.8.7]{hirsch}.

In section 4, we will also make use of the projective model
category structure $\SSets^{\Deltaop}_f$ on simplicial spaces, in
which the fibrations are given by levelwise fibrations of
simplicial sets \cite[IX 1.4]{gj}.

In \cite[7.1]{rezk}, Rezk defines a model category structure which
is obtained by localizing $\SSets^{\Deltaop}_c$ with respect to a
set of maps \cite[4.1.1]{hirsch}.  Its fibrant objects are called
Segal spaces, and they satisfy two properties: they are fibrant in
the Reedy model structure, and each Segal map $\varphi_k$ is a
weak equivalence for $k \geq 2$.  In particular, there is a
fibrant replacement functor taking any simplicial space to a Segal
space. In \cite[\S 5]{thesis}, we construct a similar functor
$L_c$ which takes a Segal precategory $X$ to a Segal space $L_cX$
which is also a Segal category such that $X_0=(L_cX)_0$.

\begin{theorem} \label{Secatc} \cite[5.1]{thesis}
There is a model category structure $\Secat_c$ on the category of
Segal precategories with the following weak equivalences,
fibrations, and cofibrations:
\begin{itemize}
\item Weak equivalences are the maps $f:X \rightarrow Y$ such that
the induced map $L_cX \rightarrow L_cY$ is a weak equivalence in
Rezk's model structure on simplicial spaces.

\item Cofibrations are the monomorphisms.  (In particular, every
Segal precategory is cofibrant.)

\item Fibrations are the maps with the right lifting property with
respect to the maps which are both cofibrations and weak
equivalences.
\end{itemize}
\end{theorem}

\begin{cor} \cite[5.13]{thesis} \label{conv}
Fibrant objects in $\Secat_c$ are Reedy fibrant Segal categories.
\end{cor}

We would like to prove the converse of this corollary.  We begin
by proving it in a more restricted setting, namely in the model
category of Segal precategories which have a fixed set $\mathcal
O$ in degree zero.

\section{Reedy fibrant Segal categories}

We begin this section by characterizing the fibrant objects in the
fixed object set case, and then we proceed to the more general
case.  We begin by briefly describing the model structure in the
more restricted situation.

We show in \cite[3.9]{simpmon} that there is a model structure
$\SSpoc$ on the category of Segal precategories which have a fixed
set $\mathcal O$ in degree zero in which the weak equivalences and
cofibrations are each defined levelwise. This model category can
be localized with respect to a set of maps to obtain a model
structure $\LSSpoc$, in which the fibrant objects are Segal
categories.

Thus, we will first prove that the fibrant objects of $\LSSpoc$
are precisely the Reedy fibrant Segal categories with the set
$\mathcal O$ in degree zero.   Note that any Reedy fibrant Segal
precategory is a fibrant object in the appropriate $\SSpoc$, and
hence that any Reedy fibrant Segal category is fibrant in the
appropriate $\LSSpoc$.  The following proposition shows that the
converse statement holds as well.

\begin{prop} \label{oc}
A fibrant object in $\SSpoc$ is fibrant in the Reedy model
structure on simplicial spaces.  In particular, a fibrant object
of $\LSSpoc$ is fibrant in the Reedy structure.
\end{prop}

\begin{proof}
Let $W$ be a fibrant object in $\SSpoc$.  We need to show that the
map $W \rightarrow \Delta[0]^t$ has the right lifting property
with respect to all levelwise acyclic cofibrations, not just the
ones in $\SSpoc$. We first consider the case where $A \rightarrow
B$ is an acyclic cofibration where $A$ and $B$ are Segal
precategories, say in $\mathcal{SS}p_{\mathcal O',c}$ for some set
$\mathcal O' \neq \mathcal O$.


Using the 0-skeleta $\sk_0(A)$ and $\sk_0(W)$, we can define a
simplicial space $A'$ as the pushout
\[ \xymatrix{\sk_0(A) \ar[r] \ar[d] & \sk_0(W) \ar[d] \\
A \ar[r] & A'.} \] Note in particular that the induced map $A'_0
\rightarrow W_0$ is an isomorphism.  Now, we can define $B'$ as a
pushout
\[ \xymatrix{A \ar[r] \ar[d] \ar[r] & A' \ar[d] \\
B \ar[r] & B'} \] and thus a diagram
\[ \xymatrix{A \ar[r] \ar[d] & A' \ar[r] \ar[d] & W \ar[d] \\
B \ar[r] \ar[drr] & B' \ar[r] & \amalg_{\mathcal O} \Delta[0]^t \ar[d] \\
&& \Delta [0]^t } \] Because it is defined as a pushout along an
acyclic cofibration in the Reedy structure, the map $A'
\rightarrow B'$ is also an acyclic cofibration, and $A'_0 \cong
B'_0 \cong \amalg_\mathcal O \Delta [0]^t$. Therefore there exists
a lift $B' \rightarrow W$, from which there exists a lift $B
\rightarrow W$.

%
%

Now, suppose that $A \rightarrow B$ is an acyclic cofibration
between simplicial spaces which are not necessarily Segal
precategories.  Since $W$ and $\Delta[0]^t$ are Segal
precategories, we can factor our diagram as
\[ \xymatrix{A \ar[r] \ar[d] & \widetilde{A} \ar[r] \ar[d] & W \ar[d] \\
B \ar[r] & \widetilde{B} \ar[r] & \Delta[0]^t} \] where
$\widetilde{A}$ and $\widetilde{B}$ are obtained from $A$ and $B$,
respectively, by collapsing the space in degree zero to its
components.  Then, we obtain a lifting from the previous argument.

Since fibrant objects in $\LSSpoc$ are fibrant in $\SSpoc$, the
second statement of the proposition follows as well.
\end{proof}

Using this result, we turn to the converse to Corollary
\ref{conv}.

%

\begin{theorem} \label{main}
Any Reedy fibrant Segal category is fibrant in $\Secat_c$.
\end{theorem}

\begin{proof}
Let $W$ be a Reedy fibrant Segal category and suppose that $f:A
\rightarrow B$ is a generating acyclic cofibration in $\Secat_c$.
We need to show that the map $W \rightarrow \Delta[0]^t$ has the
right lifting property with respect to the map $f$.  We know that
it has the right lifting property with respect to any such $f$
which preserves a fixed object set $\mathcal O$, by Proposition
\ref{oc}. Therefore, we assume that $f$ is a monomorphism but is
not surjective.

Choose $b \in B_0$ which is not in the image of $f:A \rightarrow
B$. Since $f$ is a weak equivalence, we know that $b$ is
equivalent in $L_cB$ to $f(a)$ for some $a \in (L_cA)_0=A_0$.
Define $(L_cB)_a$ to be the full Segal subcategory of $L_cB$ whose
objects are in the essential image of $a$.  Let $B_a$ be the
sub-simplicial space of $B$ whose image is in $(L_cB)_a$. Note
that $(B_a)_0=((L_cB)_a)_0$.

Now, define $A_a$ to be the Segal precategory which has as 0-space
the union of $(B_a)_0$ and $A_0$ and for which
\[ (A_a)_n(a_0, \ldots ,a_n)= A_n(a, \ldots ,a) \] for all $a_i
\in (A_a)_0$.  Letting $a$ also denote the doubly constant
simplicial space given by $a$, we can then define $A_1$ to be a
pushout given by
\[ \xymatrix{a \ar[r] \ar[d] & A_a \ar[d] \\
A \ar[r] & A_1.} \]  Notice that the map $A \rightarrow A_1$ has a
section and that we can factor $f$ as the composite $A \rightarrow
A_1 \rightarrow B$.

We now repeat this process by choosing a $b'$ which is not in the
image of the map $A_1 \rightarrow B$ and a corresponding $a'$, and
continue to do so, perhaps infinitely many times, and take a
colimit to obtain a Segal precategory $\widehat{A}$ such that the
map $f$ factors as $A \rightarrow \widehat{A} \rightarrow B$ and
there is a section $\widehat{A} \rightarrow A$.  Furthermore,
notice that $\widehat A_0 = B_0$ and that the map $\widehat A
\rightarrow B$ is object-preserving.  Notice also that it is an
acyclic cofibration.  Therefore, the dotted-arrow lift exists in
the following diagram:
\[ \xymatrix{\widehat A \ar[r] \ar[d] & W \ar[d] \\
B \ar[r] \ar@{-->}[ur] & \ast} \] which implies, using the section
$\widehat A \rightarrow A$, that there is a dotted arrow lift in
the diagram
\[ \xymatrix{A \ar[r] \ar[d] & W \ar[d] \\
B \ar[r] \ar@{-->}[ur] & \ast}. \]
\end{proof}

\begin{note}
One might wonder why, in the general case, the model category
structure $\Secat_c$ is not defined as a localization of a model
structure with object-preserving weak equivalences, as it is in
the fixed object set situation.  We prove in \cite[3.12]{thesis}
that it is impossible, simply because the more basic model
structure cannot exist. Therefore, we cannot use the tools of
localized model structures.  Much of the difficulty with working
with $\Secat_c$ arises from this fact.
\end{note}

\section{Projective fibrant Segal categories}

In this section, we show that the analogous statement holds in our
other model category structure, $\Secat_f$.  In this structure,
the weak equivalences are the same as those of $\Secat_c$, but the
fibrations and cofibrations are defined differently.  They are
technical to describe, so we will refer the interested reader to
\cite[\S 4, 7.1]{thesis} for a complete description.  The
fibrations should be thought of heuristically as those we would
obtain via a localization of a model structure given by levelwise
weak equivalences and fibrations, although, as mentioned at the
end of the previous section, such a construction is impossible in
the general case.

%
%

As before, we begin by considering the special case where we have
a fixed set $\mathcal O$ in degree zero. Analogously to the
situation in the previous section, there exists a model category
structure $\SSpof$ on simplicial spaces with a fixed set $\mathcal
O$ in degree zero in which the weak equivalences and fibrations
are levelwise. Localizing this model structure with respect to a
set of maps results in a model structure $\LSSpof$ in which the
fibrant objects are Segal categories with $\mathcal O$ in degree
zero \cite[3.8]{simpmon}.

We can prove the following result using the same techniques that
we used to prove Proposition \ref{oc}.

\begin{prop}
The fibrant objects in $\LSSpof$ are precisely the Segal
categories which have the set $\mathcal O$ in degree zero and are
fibrant in the projective model structure on simplicial spaces.
\end{prop}

As in the Reedy case, there is a corresponding functor $L_f$
(where ``Segal spaces" are now obtained as a localization in the
projective, rather than the Reedy, model structure), and in
\cite[\S 7]{thesis} we show that the weak equivalences are
actually independent of whether we define them in terms of $L_c$
or $L_f$.

\begin{theorem}
Segal categories which are fibrant in the projective model
category structure on simplicial spaces are fibrant in $\Secat_f$.
\end{theorem}

\begin{proof}
The argument given for $\Secat_c$ still holds in this case.
\end{proof}

\end{document}